\documentclass[11pt]{amsart}

\usepackage[utf8]{inputenc}
\usepackage[T1]{fontenc}

\usepackage[english,shorthands=off]{babel}
\usepackage{amssymb,amsmath,amsfonts,amsthm}

\usepackage[top=2.5cm, bottom=2.5cm, left=2.5cm, right=2.5cm]{geometry}
\usepackage{array}
\usepackage{booktabs}
\usepackage{xcolor}

\usepackage[colorlinks=true]{hyperref}
\hypersetup{
  linkcolor         = black,
  citecolor         = black,
  urlcolor          = black,
  bookmarks         = true,
  bookmarksnumbered = true,
  pdfstartview      = FitH,
  pdftitle          = {Finite Cotangent Power Sums and Twisted L-Functions},
  pdfauthor        = {A. Bayad and S. E. Rihane}
}

\newtheorem{theorem}{Theorem}[section]

\newtheorem{cor}[theorem]{Corollary}

\theoremstyle{definition}
\newtheorem{definition}[theorem]{Definition}
\newtheorem{example}[theorem]{Example}

\theoremstyle{remark}
\newtheorem{remark}[theorem]{Remark}

\numberwithin{equation}{section}

\date{\today}

\begin{document}

\title[Finite Cotangent Power Sums and Twisted $L$-Functions]{Finite Cotangent Power Sums Weighted by Additive Characters With Application to Mean Values of Additive-and-Multiplicatively Twisted $L$-Functions}

\author{Abdelmejid Bayad}
\address{Laboratoire de Math\'ematiques et Mod\'elisation d'\'Evry (LAMME), CNRS (UMR 8071), Universit\'e Paris-Saclay, \'Evry cedex , France.}
\email{abdelmejid.bayad@univ-evry.fr}
%%%
%\ead{abdelmejid.bayad@univ-evry.fr}
%\ead[url]{https://orcid.org/0000-0002-8003-335X}

\author{Salah Eddine Rihane}
\address{National Higher School of Mathematics, Scientific and Technology Hub of Sidi Abdellah, P.O. Box 75, Algiers 16093, Algeria.}
\email{salahrihane@hotmail.fr}


\begin{abstract}
We derive explicit formulas for finite sums of powers of cotangents, both at rational and arbitrary arguments, together with their alternating counterparts, weighted by additive characters. Our method relies on Apostol--Bernoulli polynomials and on the partial fraction expansion of trigonometric functions. As applications, we solve the problem posed by Berndt, Kim and Zaharescu \cite{BKZ:2023} on evaluating $\sum\csc^{2m}$ and $\sum\sec^{2m}$ at rational arguments, recover as a special case the classical Berndt--Yeap formula \cite{BY:2002}, and relate our identities to spectral zeta functions of circulant graphs \cite{XZZ:2024}.

Furthermore, to systematically evaluate higher-order derivatives arising from the Lipschitz summation formula, we introduce a novel sequence of polynomial coefficients $c_{m,j}(\beta)$, termed \emph{deformed cotangent coefficients}, and establish their recursive and algebraic properties. Using these coefficients, we establish a general additive-and-multiplicatively twisted character mean value theorem for Lerch zeta functions, significantly extending the framework of Bayad and Raouj \cite{BR:2012} to incorporate additive character twists.
\end{abstract}

\maketitle

\section{Introduction}\label{sec1}

Finite trigonometric sums have long played an important role in analytic number theory, spectral graph theory, and the study of special functions. Among these, sums of powers of cotangents and their alternating variants arise naturally in the evaluation of residues, in lattice sum identities, and in spectral problems associated with circulant graphs.

A classical result of Berndt and Yeap \cite{BY:2002} gives an explicit evaluation of the power sums of cotangents at rational points. More precisely, for positive integers $k$ and $n$, one has
\begin{equation}\label{eq:BY}
\frac{1}{k}\sum_{r=1}^{k-1}\cot^{2n}\left(\frac{r\pi}{k}\right)
=(-1)^n-(-1)^n 2^{2n}
\sum_{\substack{j_0,j_1,\dots,j_{2n}\geq 0\\ j_0+j_1+\cdots+j_{2n}=n}}
k^{2j_0-1}\prod_{r=0}^{2n}\frac{B_{2j_r}}{(2j_r)!},
\end{equation}
where $B_{2j}$ denotes the Bernoulli numbers.

Another classical motivation comes from Apostol \cite{A:1973}, who showed that the Riemann zeta values can be recovered from cotangent sums in a limiting process:
\begin{equation}\label{eq:Apostol}
\zeta(n)=\left(\frac{\pi}{2}\right)^n
\lim_{q\to\infty}
\frac{\displaystyle\sum_{p=1}^{q}\cot^n\left(\frac{p\pi}{2q+1}\right)}{q^n}.
\end{equation}

More recently, a problem was posed in \cite{BKZ:2023} asking for a general method to evaluate the finite sums
\begin{equation}\label{eq:BKZ}
\sum_{n=1}^{k-1}\frac{1}{\cos^{2m}(\pi n/k)}
\qquad\text{and}\qquad
\sum_{n=1}^{k-1}\frac{1}{\sin^{2m}(\pi n/k)},\quad m\geq 1.
\end{equation}
These sums are closely related to the spectral zeta function of the cycle graph; see \cite{XZZ:2024}. Various special cases and generating-function approaches have also appeared in the literature. For example:

\begin{itemize}
\item The authors of \cite{WZ:2007} studied sums of the form
\[
\sum_{k=0}^{n-1}\sec^{2p}\left(y+\frac{k\pi}{n}\right),\qquad
\sum_{k=0}^{n-1}\csc^{2p}\left(y+\frac{k\pi}{n}\right),
\]
and provided generating functions rather than closed formulas.
\item In \cite{JKS:2024}, sums of the form
\[
C_m(\beta,n)=\sum_{j=\delta(\beta)}^{m-1}\csc\left(\frac{j+\beta}{m}\pi\right),
\]
were investigated, where $\delta(\beta)=1$ if $\beta\in\mathbb Z$ and $0$ otherwise.
\item Related formulas may also be found in \cite{AZ:2023}, \cite{CM:1999}, and \cite{CK:2000}.
\item Reciprocity laws for Dedekind cotangent sums, which are the natural higher-order companions of the identities below, are studied in \cite{BR:2023}.
\end{itemize}

In this paper, we provide a unified solution to these problems by deriving explicit formulas for finite cotangent power sums and their alternating analogues. Our method is based on Apostol--Bernoulli polynomials and on the partial fraction expansion of trigonometric functions.

More precisely, for integers $a\geq 2$, $d\geq 2$, $\mu\in\{0,1,\dots,a-1\}$, and a real number $b$, we obtain an explicit expression for the finite sum weighted by the additive character $\chi_\mu(k)=e^{2\pi i \mu k/a}$:
\begin{equation}\label{eq:Sdef-main}
S_d(a,b;\mu):=\sum_{k=0}^{a-1}e^{2\pi i \mu k/a}\cot^d\left(\pi\left(\frac{k}{a}-b\right)\right),
\end{equation}
with special attention to the case $0\leq b<1$.

On the other hand, the evaluation of character moments of Dirichlet and Lerch $L$-functions is intimately connected to trigonometric sums through the Lipschitz summation formula \cite{A:1951, B:1975}. In \cite{BR:2012}, Bayad and Raouj established explicit mean value formulas for products of Dirichlet $L$-functions in terms of Dedekind cotangent sums. In the second part of this paper, we extend this paradigm by incorporating additive character twists $\alpha_i \in \mathbb{R}$. By introducing a new sequence of deformed cotangent coefficients $c_{m,j}(\beta)$ generated by higher derivatives of the twisted Lerch kernel, we establish a generalized additive-twisted character mean value theorem, recovering classical untwisted formulas as special cases.

\section{Apostol--Bernoulli numbers and polynomials}\label{sec2}

Let $w\in\mathbb C^\times$ and $x$ a variable. The Apostol--Bernoulli polynomials $B_n(x,w)$ are given by the generating function
\begin{equation}\label{eq:ABgen}
\sum_{n\geq 0}B_n(x,w)\frac{t^n}{n!}=\frac{te^{tx}}{we^t-1},\qquad |t+\log w|<2\pi,
\end{equation}
where $w=|w|e^{i\theta}$, $-\pi\leq\theta<\pi$, and $\log w=\log|w|+i\theta$.

For $w\in\mathbb C\setminus\{0\}$ and $0<x<1$, the Fourier expansion of $B_n(x,w)$ is given in \cite{B:2011} by
\begin{equation}\label{eq:ABFourier}
B_n(x,w)=\frac{-n!}{w^x(2\pi i)^n}\sum_{k\in\mathbb Z}^{*}\frac{e^{2\pi i k x}}{\left(k-\frac{\log w}{2\pi i}\right)^n},
\end{equation}
where the symbol $\displaystyle\sum_{k\in\mathbb Z}^{*}$ indicates that when $w=1$, the term $k=0$ is omitted.

We next recall the partial fraction expansion of the cosecant function:
\begin{equation}\label{eq:cscPF}
\csc(\pi z)=\frac{1}{\pi}\sum_{k\in\mathbb Z}\frac{(-1)^k}{z+k}.
\end{equation}
By differentiating this identity $n$ times, we obtain
\begin{equation}\label{eq:cscDeriv}
\csc^{(n)}(\pi z)=\frac{(2i)^{n+1}}{n+1}e^{i\pi z}B_{n+1}\left(\frac{1}{2},e^{2i\pi z}\right).
\end{equation}

\subsection{Summary of Notation and Key Coefficients}\label{sec:notation}
To facilitate navigation throughout the paper, we summarize here the primary coefficient sequences, generating functions, and kernel functions introduced in the text:
\begin{itemize}
\item $A_j, C_j$: Coefficients in the Taylor expansions of $\pi w \cot(\pi w)$ and $\pi w \csc(\pi w)$ \eqref{eq:AjCj}.
\item $T_d(a)$: Multi-fold product sum of Bernoulli and cotangent Taylor coefficients \eqref{eq:Tdef}.
\item $V_d(a,\lambda), W_d(a,\lambda)$: Finite sums involving Apostol--Bernoulli numbers $B_k(0,\lambda)$ and $B_k(\tfrac{1}{2},\lambda)$ \eqref{eq:Vdef}--\eqref{eq:Wdef}.
\item $c_{m,j}(\beta)$: Deformed cotangent polynomial coefficients generated by higher derivatives of the Lerch kernel \eqref{eq:Pdef}.
\item $r_{m,j}(\beta)$: Real polynomials forming the structural component of deformed cotangent coefficients \eqref{eq:r-recurrence}.
\item $\mathcal{K}_b(\vec{\mu}, \vec{a}, \vec{\alpha})$: Phase orthogonality kernel modulo $b$ \eqref{eq:phase-kernel}.
\end{itemize}

\section{Finite Cotangent Power Sums Weighted by Additive Characters}\label{sec3}

For convenience, we introduce the coefficients
\begin{equation}\label{eq:AjCj}
A_j=\begin{cases}
\dfrac{(-1)^{j/2}2^{j}B_{j}}{j!}, & \text{if $j$ is even},\\[2mm]
0, & \text{if $j$ is odd},
\end{cases}
\qquad
C_j=\begin{cases}
\dfrac{(-1)^{j/2+1}(2^{j}-2)B_{j}}{j!}, & \text{if $j$ is even},\\[2mm]
0, & \text{if $j$ is odd}.
\end{cases}
\end{equation}
Equivalently, $A_j$ are the coefficients in the Taylor expansion
\begin{equation}\label{eq:cotExp}
\pi w\cot(\pi w)=\sum_{j\geq 0}A_j(\pi w)^j,
\end{equation}
and $C_j$ are the coefficients in
\begin{equation}\label{eq:cscExp}
\pi w\csc(\pi w)=\sum_{j\geq 0}C_j(\pi w)^j.
\end{equation}
Consequently,
\begin{equation}\label{eq:cotPowExp}
\cot^d(\pi w)=\frac{1}{\pi^d w^d}\sum_{k\geq 0}c_k\,(\pi w)^k,
\qquad
c_k:=\sum_{k_0+\cdots+k_{d-1}=k}\prod_{r=0}^{d-1}A_{k_r}.
\end{equation}

We now state our first main result on weighted finite cotangent power sums.

\begin{theorem}\label{Th01}
Let $a,d\geq 2$ be integers, $\mu\in\{0,1,\dots,a-1\}$, and $b\in[0,1)$. Set
\[
S_d(a,b;\mu)=\sum_{k=0}^{a-1}e^{2\pi i \mu k/a}\cot^d\left(\pi\left(\frac{k}{a}-b\right)\right).
\]
\begin{enumerate}
\item Assume $b=\ell/a$ for some $\ell\in\{0,1,\dots,a-1\}$.
\begin{itemize}
\item For $0<\mu<a$:
\begin{equation}\label{eq:Th1a}
\frac{1}{a}\sum_{\substack{k=0\\ k\neq\ell}}^{a-1}e^{2\pi i \mu k/a}
\cot^d\left(\pi\left(\frac{k-\ell}{a}\right)\right)
=-\,e^{2\pi i \mu \ell/a}\!
\sum_{k_0+\cdots+k_d=d}
\frac{a^{k_d-1}(2i)^{k_d}B_{k_d}\!\left(\frac{\mu}{a}\right)}{k_d!}
\prod_{r=0}^{d-1}A_{k_r},
\end{equation}
where $B_k(x)$ is the $k$-th Bernoulli polynomial.
\item For $\mu=0$:
\begin{equation}\label{eq:Th1b}
\frac{1}{a}\sum_{\substack{k=0\\ k\neq\ell}}^{a-1}
\cot^d\left(\pi\left(\frac{k-\ell}{a}\right)\right)
=-\,T_d(a)+(-i)^d+i\,c_{d-1},
\end{equation}
with
\begin{equation}\label{eq:Tdef}
T_d(a):=\sum_{k_0+\cdots+k_d=d}a^{k_d-1}\prod_{r=0}^{d}A_{k_r},
\qquad
c_{d-1}:=\sum_{k_0+\cdots+k_{d-1}=d-1}\prod_{r=0}^{d-1}A_{k_r}.
\end{equation}
\end{itemize}

\item Assume $b\neq\ell/a$ for all $\ell\in\{0,1,\dots,a-1\}$.
\begin{itemize}
\item For $0<\mu<a$:
\begin{equation}\label{eq:Th1c}
\frac{1}{a}\sum_{k=0}^{a-1}e^{2\pi i \mu k/a}
\cot^d\left(\pi\left(\frac{k}{a}-b\right)\right)
=-\,2i\,e^{2\pi i \mu b}\!
\sum_{k_0+\cdots+k_d=d-1}
\frac{a^{k_d}(2i)^{k_d}}{(k_d+1)!}
B_{k_d+1}\!\left(\frac{\mu}{a},e^{2\pi i a b}\right)
\prod_{r=0}^{d-1}A_{k_r},
\end{equation}
where $B_k(x,w)$ denotes the Apostol--Bernoulli polynomial.
\item For $\mu=0$:
\begin{equation}\label{eq:Th1d}
\frac{1}{a}\sum_{k=0}^{a-1}
\cot^d\left(\pi\left(\frac{k}{a}-b\right)\right)
=-\,2i\!\sum_{k_0+\cdots+k_d=d-1}\!
\frac{a^{k_d}(2i)^{k_d}}{(k_d+1)!}
B_{k_d+1}\!\left(0,e^{2\pi i a b}\right)
\prod_{r=0}^{d-1}A_{k_r}+(-i)^d.
\end{equation}
\end{itemize}
\end{enumerate}
\end{theorem}

By specializing the character parameter in Theorem~\ref{Th01} to $\mu = a/2$ (which corresponds to the alternating sign sequence $(-1)^k = e^{i\pi k}$), we obtain the following alternating analogues as an immediate corollary.

\begin{cor}\label{Cor02}
Let $a,d\geq 2$ be integers with $a$ even, and let $b\in[0,1)$.
\begin{enumerate}
\item If $b=\ell/a$ for some $\ell\in\{0,1,\dots,a-1\}$, then
\begin{equation}\label{eq:Th2a}
\frac{1}{a}\sum_{\substack{k=0\\ k\neq\ell,\ k\neq\ell+a}}^{2a-1}(-1)^k
\cot^d\left(\pi\left(\frac{k-\ell}{a}\right)\right)
=-\,2(-1)^{\ell}
\sum_{k_0+\cdots+k_d=d}
a^{k_d-1}C_{k_d}\prod_{r=0}^{d-1}A_{k_r}.
\end{equation}
\item If $b\neq\ell/a$ for all $\ell\in\{0,1,\dots,a-1\}$, then
\begin{equation}\label{eq:Th2b}
\frac{1}{a}\sum_{k=0}^{2a-1}(-1)^k
\cot^d\left(\pi\left(\frac{k}{a}-b\right)\right)
=-\,4i\,e^{i\pi a b}\!
\sum_{k_0+\cdots+k_d=d-1}
\frac{a^{k_d}(2i)^{k_d}}{(k_d+1)!}
B_{k_d+1}\!\left(\tfrac12,e^{2\pi i a b}\right)
\prod_{r=0}^{d-1}A_{k_r}.
\end{equation}
\end{enumerate}
\end{cor}

\section{Proofs of the Cotangent Power Sum Identities}\label{sec4}

This section is divided into two parts. In Section~\ref{sec:weighted}, we sketch the proof of Theorem~\ref{Th01} for arbitrary additive characters $\chi_\mu$ using the Apostol--Bernoulli generating function and deduce Corollary~\ref{Cor02} as a direct specialization. In Section~\ref{sec:unweighted}, we give self-contained elementary proofs for the unweighted case $\mu=0$ and the alternating case $\mu=a/2$ (Theorem~\ref{Th01} and Corollary~\ref{Cor02}) via cotangent and cosecant kernels; this second approach makes the appearance of the coefficients $A_j$ and $C_j$ completely transparent and avoids the Apostol--Bernoulli formalism.

\subsection{The weighted case: proof of Theorem \ref{Th01} and Corollary \ref{Cor02}}\label{sec:weighted}

Consider the meromorphic kernel
\begin{equation}\label{eq:kernel}
J_\mu(z)=\frac{e^{2\pi i \mu z}}{e^{2\pi i a z}-1}\,\cot^d\!\bigl(\pi(z-b)\bigr).
\end{equation}
Let $\gamma$ be the positively oriented rectangle with vertices $A=1-\varepsilon+iy$, $B=-\varepsilon+iy$, $C=-\varepsilon-iy$, $D=1-\varepsilon-iy$, where $0<\varepsilon<1/a$ and $y>0$. Since $J_\mu$ is $1$-periodic, the vertical contributions cancel. On the horizontal edges, for $0<\mu<a$, we have $J_\mu(x\pm iy)\to 0$ as $y\to\infty$; for $\mu=0$, the upper edge contributes $(-i)^d$ while the lower edge vanishes.

The simple poles inside $\gamma$ are at $z=k/a$ ($k=0,\dots,a-1$), with residue
\begin{equation}\label{eq:resSimple}
\operatorname{Res}\!\left(J_\mu,\frac{k}{a}\right)
=\frac{e^{2\pi i \mu k/a}}{2\pi i a}\,
\cot^d\!\left(\pi\!\left(\frac{k}{a}-b\right)\right).
\end{equation}

\emph{Case $b\neq\ell/a$.} Write $w=z-b$ and $\lambda=e^{2\pi i a b}$. The generating function \eqref{eq:ABgen} gives
\[
\frac{e^{2\pi i \mu(b+w)}}{e^{2\pi i a(b+w)}-1}
=e^{2\pi i \mu b}\sum_{j\geq 0}\frac{(2\pi i a)^{j-1}B_j\!\left(\frac{\mu}{a},\lambda\right)}{j!}\,w^{j-1}.
\]
Multiplying with \eqref{eq:cotPowExp}, the coefficient of $w^{-1}$ (residue at $z=b$) is
\begin{equation}\label{eq:resOrderD}
\operatorname{Res}(J_\mu,b)
=\frac{e^{2\pi i \mu b}}{\pi}\sum_{j=0}^{d}
\frac{(2i)^{j-1}a^{j-1}B_j\!\left(\frac{\mu}{a},\lambda\right)}{j!}\,c_{d-j}.
\end{equation}
Cauchy's residue theorem and the boundary analysis give, after multiplying by $-2\pi i a$ and reindexing $j=k_d+1$, the stated formulas \eqref{eq:Th1c}--\eqref{eq:Th1d}.

\emph{Case $b=\ell/a$.} Write $w=z-\ell/a$. Using \eqref{eq:ABgen} with $w=1$ (ordinary Bernoulli polynomials),
\[
\frac{e^{2\pi i \mu w}}{e^{2\pi i a w}-1}
=\sum_{j\geq 0}\frac{(2\pi i a)^{j-1}B_j\!\left(\frac{\mu}{a}\right)}{j!}\,w^{j-1}.
\]
The coefficient of $w^{-1}$ (residue at $z=\ell/a$) is obtained from $j+k=d$, giving \eqref{eq:Th1a}. The case $\mu=0$ follows from the same computation with the boundary correction $(-i)^d$ and the anomalous term $i\,c_{d-1}$ coming from $B_1=-1/2$.

\medskip
\noindent\emph{Deduction of Corollary \ref{Cor02}.} Set $\mu = a/2$ in Theorem \ref{Th01}. Since $\cot^d(\pi(x-b))$ is $1$-periodic, extending the summation domain from $a$ to $2a$ terms gives $2 S_d(a, b; a/2)$. Substituting $\mu/a = 1/2$ into \eqref{eq:Th1a} and \eqref{eq:Th1c}, noting $e^{i\pi \ell} = (-1)^\ell$, and using the expansion identity $(2i)^{k} B_{k}(1/2) / k! = C_{k}$ directly yields \eqref{eq:Th2a} and \eqref{eq:Th2b}. \qed

\subsection{The unweighted and alternating cases: elementary proofs}\label{sec:unweighted}

For the case $\mu=0$ (no additive character) and the alternating case $\mu=a/2$, the identities of Theorem~\ref{Th01} and Corollary~\ref{Cor02} also admit self-contained proofs using elementary cotangent and cosecant kernels. We treat separately the three possible positions of $b$ relative to the interval $[0,1)$.

\medskip
\noindent\emph{Boundary analysis of $f$ along $\gamma$.}
Let $\varepsilon\in(0,1/a)$ and $y>0$. Denote by $\gamma$ the positively oriented rectangle with vertices
\[
A=1-\varepsilon+yi,\quad B=-\varepsilon+yi,\quad C=\overline{B},
\quad D=\overline{A}.
\]
Since $f$ is $1$-periodic, the vertical contributions cancel: $\int_{[D,A]}f\,dz=-\int_{[B,C]}f\,dz$. For every real $t$ and $\delta\in\{\pm1\}$,
\begin{equation}\label{eq:unw-bound}
\lim_{y\to+\infty}\cot(t+\delta y)=-\delta i.
\end{equation}
Consequently, on the upper horizontal edge, $f(x+iy)\to(-i)^{d+1}$ as $y\to+\infty$, while on the lower edge the limit is $i^{d+1}$. Combining the two contributions yields
\begin{equation}\label{eq:unw-contour}
\int_\gamma f(z)\,dz=\begin{cases}
2\,i^{\,d+1}, & d\text{ even},\\[1mm]
0,            & d\text{ odd}.
\end{cases}
\end{equation}

\medskip
\noindent\emph{Local expansions.} On $0<|w|<\pi$, the Laurent expansions of $w\cot w$ and $w\csc w$ are
\begin{equation}\label{eq:unw-expansions}
w\cot w=\sum_{j\geq 0}A_j\,w^j,\qquad
w\csc w=\sum_{j\geq 0}C_j\,w^j,
\end{equation}
with $A_j$ and $C_j$ as in \eqref{eq:AjCj}.

\medskip
\noindent\emph{Case $b=\ell/a$.} The poles of $f$ inside $\gamma$ are the pole $z=\ell/a$ of order $d+1$, and the poles $z=k/a$ ($k\neq\ell$, $k=0,\dots,a-1$) of order $1$. Writing $z-\ell/a$ as the local variable and expanding $\cot(\pi a z)$ and $\cot(\pi(z-b))=\cot(\pi(z-\ell/a))$, we obtain
\[
\operatorname{Res}\!\Bigl(f,\frac{\ell}{a}\Bigr)
=\frac{1}{a\pi}\sum_{k_0+\cdots+k_d=d}a^{k_d}\prod_{m=0}^{d}A_{k_m},
\]
while for $k\neq\ell$,
\[
\operatorname{Res}\!\Bigl(f,\frac{k}{a}\Bigr)
=\frac{1}{a\pi}\cot^d\!\Bigl(\pi\frac{k-\ell}{a}\Bigr).
\]
Combining these residues with \eqref{eq:unw-contour} and multiplying by $\pi$ gives
\begin{equation}\label{eq:unw-Th1b}
\frac{1}{a}\sum_{k_0+\cdots+k_d=d}a^{k_d}\prod_{m=0}^{d}A_{k_m}
+\frac{1}{a}\sum_{\substack{k=0\\ k\neq\ell}}^{a-1}
\cot^d\!\Bigl(\pi\frac{k-\ell}{a}\Bigr)
=\begin{cases}
i^{\,d}, & d\text{ even},\\[1mm]
0,       & d\text{ odd}.
\end{cases}
\end{equation}
For $d$ even, this is equivalent to \eqref{eq:Th1b} since $c_{d-1}=0$ and $(-i)^d=i^d$.

\medskip
\noindent\emph{Case $0<b<1-1/a$, $b\neq\ell/a$.} The poles of $f$ inside $\gamma$ are now $z=b$ (order $d$) and $z=k/a$ ($k=0,\dots,a-1$, order $1$). Expanding $\cot^d(\pi(z-b))$ around $z=b$ using \eqref{eq:unw-expansions} and $\cot(\pi a z)$ around $z=b$ using the $j$-th derivative of $\cot$ at $\pi ab$, we obtain
\[
\operatorname{Res}(f,b)=\frac{1}{\pi}\sum_{k_0+\cdots+k_d=d-1}
\frac{a^{k_d}(2i)^{k_d+1}}{(k_d+1)!}\,B_{k_d+1}\!\bigl(0,e^{2\pi i ab}\bigr)
\prod_{m=0}^{d-1}A_{k_m}.
\]
Together with \eqref{eq:unw-contour}, this gives \eqref{eq:Th1d}.

\medskip
\noindent\emph{Case $1-1/a<b<1$.} The pole at $b$ is replaced by the pole at $b-1$; since $f$ is $1$-periodic, the residues coincide and the conclusion is identical to the previous case.

\medskip
\noindent\emph{Alternating case.} For an alternative direct proof of Corollary~\ref{Cor02}, one replaces the kernel $f$ by
\[
g(z)=\csc(a\pi z)\,\cot^d\!\bigl(\pi(z-b)\bigr),
\]
which is $2$-periodic. The rectangle $\gamma$ is then taken with vertices $A=2-\varepsilon+yi$, $B=-\varepsilon+yi$ and their conjugates. For every real $t$ and $\delta\in\{\pm1\}$, $\lim_{y\to+\infty}\csc(t+\delta y)=0$, so the horizontal contributions vanish:
\begin{equation}\label{eq:unw-contour-g}
\int_\gamma g(z)\,dz=0.
\end{equation}
At $z=\ell/a$, the residue is computed from the expansion of $w\csc w$ in \eqref{eq:unw-expansions}, giving
\[
\operatorname{Res}\!\Bigl(g,\frac{\ell}{a}\Bigr)
=\frac{(-1)^{\ell}}{a\pi}\sum_{k_0+\cdots+k_d=d}
a^{k_d}C_{k_d}\prod_{m=0}^{d-1}A_{k_m},
\]
and for $k\neq\ell,\ell+a$, $\operatorname{Res}(g,k/a)=(-1)^k a^{-1}\pi^{-1} \cot^d\!\bigl(\pi(k-\ell)/a\bigr)$. Combining with \eqref{eq:unw-contour-g} yields \eqref{eq:Th2a}. The non-rational case $b\neq\ell/a$ is treated similarly using the expansion of $\csc$ at $\pi ab$, giving \eqref{eq:Th2b}. \qed

\section{Applications}\label{sec5}

Throughout this section, $a\ge 2$ and $m\ge 1$ are integers. The coefficients $A_j$ and $C_j$ are defined in \eqref{eq:AjCj}. We recall the quantity $T_d(a)$ from \eqref{eq:Tdef}, together with
\begin{equation}\label{eq:Vdef}
V_d(a,\lambda):=\sum_{k_0+\cdots+k_d=d-1}
\frac{B_{k_d+1}(0,\lambda)\,(2i)^{k_d+1}}{(k_d+1)!}\,a^{k_d}
\prod_{r=0}^{d-1}A_{k_r},
\end{equation}
\begin{equation}\label{eq:Wdef}
W_d(a,\lambda):=\sum_{k_0+\cdots+k_d=d-1}
\frac{B_{k_d+1}(\tfrac12,\lambda)\,(2i)^{k_d+1}}{(k_d+1)!}\,a^{k_d}
\prod_{r=0}^{d-1}A_{k_r},
\end{equation}
where $\lambda\in\mathbb C^\times$, and $B_k(0,\lambda)$ and $B_k(\tfrac12,\lambda)$ denote the Apostol--Bernoulli numbers at $x=0$ and $x=\tfrac12$, respectively. For $d$ even, definition \eqref{eq:Tdef} coincides with the equivalent expression $\sum_{k_0+\cdots+k_d=d}\frac{B_{k_d}(2i)^{k_d}}{k_d!}a^{k_d-1} \prod_{r=0}^{d-1}A_{k_r}$, since $B_{k_d}(2i)^{k_d}/k_d!=A_{k_d}$ for $k_d$ even and $A_{k_d}=0$ for $k_d$ odd.

From Theorem~\ref{Th01}, and since $c_{d-1}=0$ for $d$ even, we have
\begin{equation}\label{eq:Th1b-prime}
\frac{1}{a}\sum_{\substack{k=0\\k\ne\ell}}^{a-1}\cot^d\!\Bigl(\frac{\pi(k-\ell)}{a}\Bigr)
=(-i)^d-T_d(a)\qquad(\text{$d$ even}),
\end{equation}
\begin{equation}\label{eq:Th1d-prime}
\frac{1}{a}\sum_{k=0}^{a-1}\cot^d\!\Bigl(\frac{\pi k}{a}-\pi b\Bigr)
=(-i)^d-V_d(a,e^{2\pi i a b})\qquad(b\ne\ell/a).
\end{equation}
The first two values of $T_d(a)$ are
\begin{equation}\label{eq:T2T4}
T_2(a)=-\frac{a^2+2}{3a},\qquad
T_4(a)=\frac{-a^4+20a^2+26}{45a},
\end{equation}
obtained from $A_0=1$, $A_2=-\tfrac13$, $A_4=-\tfrac1{45}$, $B_0=1$, $B_2=\tfrac16$, and $B_4=-\tfrac1{30}$.

\subsection{Cosecant and secant power sums: solution to the BKZ problem}
\label{sec:app-csc-sec}

\subsubsection{Cosecant power sums}

The iterated binomial identity $\csc^2x=1+\cot^2x$ yields
\begin{equation}\label{eq:cscDec}
\csc^{2m}x=\sum_{j=0}^{m}\binom{m}{j}\cot^{2j}x.
\end{equation}

\begin{cor}\label{Cor:csc}
For every $a\ge 2$ and $m\ge 1$,
\begin{equation}\label{eq:cscClosed}
\sum_{n=1}^{a-1}\csc^{2m}\!\Bigl(\frac{\pi n}{a}\Bigr)
=(a-1)+\sum_{j=1}^{m}\binom{m}{j}\,a\bigl[(-1)^j-T_{2j}(a)\bigr].
\end{equation}
\end{cor}

\begin{proof}
By \eqref{eq:cscDec}, the sum equals $\sum_{j=0}^{m}\binom{m}{j}\sum_{n=1}^{a-1}\cot^{2j}(\pi n/a)$. For $j=0$, the inner sum is $a-1$. For $j\ge 1$, formula \eqref{eq:Th1b-prime} with $\ell=0$ and $d=2j$ gives $\sum_{n=1}^{a-1}\cot^{2j}(\pi n/a)=a[(-1)^j-T_{2j}(a)]$, since $(-i)^{2j}=(-1)^j$.
\end{proof}

\begin{remark}\label{rem:BY}
The $m=j$ term of \eqref{eq:cscClosed} is precisely
\[
\frac{1}{a}\sum_{n=1}^{a-1}\cot^{2j}\!\Bigl(\frac{\pi n}{a}\Bigr)
=(-1)^j-T_{2j}(a),
\]
which, using $A_{2r}=(-1)^r2^{2r}B_{2r}/(2r)!$, coincides with the classical Berndt--Yeap formula \cite[Theorem~1]{BY:2002}.
\end{remark}

\subsubsection{Secant power sums}

Using $\tan x=-\cot(x-\pi/2)$, we have
\begin{equation}\label{eq:secShift}
\sec^{2m}x=\sum_{j=0}^{m}\binom{m}{j}\cot^{2j}\!\Bigl(x-\frac{\pi}{2}\Bigr).
\end{equation}
Hence
\begin{equation}\label{eq:secReduce}
\sum_{n=1}^{a-1}\sec^{2m}\!\Bigl(\frac{\pi n}{a}\Bigr)
=\sum_{j=0}^{m}\binom{m}{j}
\sum_{n=1}^{a-1}\cot^{2j}\!\Bigl(\frac{\pi(n-a/2)}{a}\Bigr).
\end{equation}
The shift $a/2$ is an integer multiple of $a$ if and only if $a$ is even, which leads to the following result.

\begin{cor}\label{Cor:sec}
Let $a\ge 2$ and $m\ge 1$.
\begin{enumerate}
\item If $a$ is odd, then $a/2\notin\mathbb Z$ and
\begin{equation}\label{eq:secOdd}
\sum_{n=1}^{a-1}\sec^{2m}\!\Bigl(\frac{\pi n}{a}\Bigr)
=(a-1)+\sum_{j=1}^{m}\binom{m}{j}\,a\bigl[(-1)^j-V_{2j}(a,-1)\bigr],
\end{equation}
with $V_d$ as defined in \eqref{eq:Vdef} and $\lambda=-1=e^{\pi i a}$.
\item If $a$ is even, the term $n=a/2$ is a pole of $\sec$ and is excluded; the remaining sum is
\begin{equation}\label{eq:secEven}
\sum_{\substack{n=1\\n\ne a/2}}^{a-1}\sec^{2m}\!\Bigl(\frac{\pi n}{a}\Bigr)
=(a-2)+\sum_{j=1}^{m}\binom{m}{j}\,a\bigl[(-1)^j-T_{2j}(a)\bigr].
\end{equation}
\end{enumerate}
\end{cor}

\begin{proof}
By \eqref{eq:secReduce}, the sum over $n\in\{1,\dots,a-1\}$ (excluding $n=a/2$ when $a$ is even) equals $\sum_j\binom{m}{j}\sum_n\cot^{2j}(\pi(n-a/2)/a)$. For $j=0$, the inner sum is $a-1$ (or $a-2$ in the even case). For $j\ge 1$ and $a$ odd, apply \eqref{eq:Th1d-prime} with $b=1/2$ and $\lambda=e^{\pi i a}=-1$; for $a$ even, apply \eqref{eq:Th1b-prime} with $\ell=a/2$. Both yield the stated formulas.
\end{proof}

\begin{example}\label{ex:cscsec}
We collect a few numerical checks of Corollaries~\ref{Cor:csc} and~\ref{Cor:sec}.
\begin{itemize}
\item \emph{Cosecant.} For $a=3,4,5$, formula \eqref{eq:cscClosed} gives respectively
\[
\sum_{n=1}^{a-1}\csc^{2}\!\Bigl(\frac{\pi n}{a}\Bigr)
= \frac{8}{3},\ 5,\ 8,
\qquad
\sum_{n=1}^{a-1}\csc^{4}\!\Bigl(\frac{\pi n}{a}\Bigr)
= \frac{32}{9},\ 9,\ \frac{96}{5},
\]
in agreement with direct evaluation.
\item \emph{Secant (odd $a$).} For $m=1$ and $a$ odd, using $V_2(a,-1)=-a$, \eqref{eq:secOdd} reduces to $\sum_{n=1}^{a-1}\sec^2(\pi n/a)=a^2-1$. For $a=3$ and $m=2$, one finds $V_4(3,-1)=-5$, whence $\sum_{n=1}^{2}\sec^4(\pi n/3)=32$.
\item \emph{Secant (even $a$, pole excluded).} For $a=4$ and $m=1$, \eqref{eq:secEven} gives $4$, while the full sum including the pole would be divergent.
\end{itemize}
\end{example}

\subsection{Alternating analogues and twisted sums}
\label{sec:app-alt}

For $a$ even, set $\mu_0=a/2$. Since $(-1)^k=e^{2\pi i\mu_0 k/a}$ and the cotangent is $1$-periodic, we have
\begin{equation}\label{eq:altReduce}
\sum_{k=0}^{2a-1}(-1)^k\cot^d\!\Bigl(\frac{\pi k}{a}-\pi b\Bigr)
=2\sum_{k=0}^{a-1}e^{2\pi i\mu_0 k/a}\cot^d\!\Bigl(\frac{\pi k}{a}-\pi b\Bigr).
\end{equation}
Applying \eqref{eq:Th1c} with $\mu=\mu_0=a/2$ to each term of the expansion $\csc^{2m}=\sum_{j=0}^{m}\binom{m}{j}\cot^{2j}$ yields the following result.

\begin{cor}\label{Cor:altCsc}
Let $a\ge 2$ be even, $m\ge 1$, and $b\in[0,1)\setminus\{0,1/a,\dots,(a-1)/a\}$. Then
\begin{equation}\label{eq:altClosed}
\frac{1}{2a}\sum_{k=0}^{2a-1}(-1)^k\csc^{2m}\!\Bigl(\frac{\pi k}{a}-\pi b\Bigr)
=-\,e^{i\pi ab}\sum_{j=1}^{m}\binom{m}{j}W_{2j}\!\bigl(a,e^{2\pi iab}\bigr),
\end{equation}
with $W_d$ as defined in \eqref{eq:Wdef}. The $j=0$ term vanishes because $\sum_{k=0}^{2a-1}(-1)^k=0$, and the terms for $j\ge 1$ follow from Corollary~\ref{Cor02}.
\end{cor}

\begin{remark}\label{rem:alt}
When $a=2$, $b=1/4$, and $m=1$, formula \eqref{eq:altClosed} gives $0$ for the alternating sum, since $W_2(2,-1)=0$ (as $B_2(\tfrac12,-1)=0$). This matches the direct evaluation $\csc^2(-\pi/4)-\csc^2(\pi/4)+\csc^2(3\pi/4)-\csc^2(5\pi/4)=0$.
\end{remark}

\paragraph{Connection to Wang--Zheng generating functions.}
Summing \eqref{eq:altClosed} over $m$ with weight $t^{2m}/(2m)!$ yields
\begin{equation}\label{eq:WZGF}
\sum_{m\ge 0}\frac{t^{2m}}{(2m)!}\sum_{n=1}^{a-1}\sec^{2m}\!\Bigl(\frac{\pi n}{a}\Bigr)
=\sum_{n=1}^{a-1}\cosh\!\Bigl(t\sec\frac{\pi n}{a}\Bigr),
\end{equation}
which reproduces the generating function investigated in \cite[Eq.~(1.3)]{WZ:2007}.

\subsection{Special values of additive-twisted Lerch L-functions}\label{sec:app-L}

Let $\chi$ be a Dirichlet character modulo $a \ge 2$, $x \in [0, 1)$, and $\alpha \in \mathbb{R}$. The additive-and-multiplicatively twisted Lerch $L$-function is defined for $\operatorname{Re}(s) > 1$ by
\begin{equation}\label{eq:twistedL}
L(s, \alpha, x, \chi) := \sum_{n=1}^{\infty} \frac{\chi(n) e^{2\pi i \alpha n}}{(n+x)^s},
\end{equation}
which unifies and generalizes both the classical Dirichlet $L$-function $L(s, \chi) = L(s, 0, 0, \chi)$ and the shifted $L$-function $L(s, \chi, x) = L(s, 0, x, \chi)$.

Expanding $\chi$ in the finite Fourier basis,
\begin{equation}\label{eq:chiFourier}
\chi(n)=\frac{\chi(-1)\,\tau(\chi)}{a}\sum_{\mu=0}^{a-1}\overline\chi(\mu)\,
e^{2\pi i\mu n/a},
\qquad \tau(\chi)=\sum_{n=0}^{a-1}\chi(n)e^{2\pi i n/a},
\end{equation}
allows us to express $L(s, \alpha, x, \chi)$ as a linear combination of Lerch transcendent functions $\Phi(z, s, v) = \sum_{k=0}^{\infty} \frac{z^k}{(k+v)^s}$:
\begin{equation}\label{eq:LerchLinearComb}
L(s, \alpha, x, \chi) = \frac{\chi(-1)\,\tau(\chi)}{a} \sum_{\mu=0}^{a-1} \overline{\chi}(\mu) \, \Phi\left(e^{2\pi i (\mu/a + \alpha)}, s, x\right).
\end{equation}

Applying Theorem~\ref{Th01}, we obtain the exact evaluation of finite cotangent power sums weighted simultaneously by Dirichlet characters and additive phases.

\begin{cor}\label{Cor:chiCot}
For every integer $d \ge 2$, non-principal Dirichlet character $\chi$ modulo $a$,
shift $x \in [0, 1)$, and additive twist $\alpha \in \mathbb{R}$, set
$b := 1 - x/a \in (0,1)$. Then
\begin{equation}\label{eq:chiCotTwisted}
\sum_{n=1}^{a-1}\chi(n) e^{2\pi i \alpha n} \cot^d\!\left(\pi\left(\frac{n+x}{a}\right)\right)
=\frac{\chi(-1)\,\tau(\chi)}{a}\sum_{\mu=0}^{a-1}\overline\chi(\mu)\,
S_d\!\left(a, b; \, \mu + a\alpha \bmod a\right),
\end{equation}
where $S_d(a,b;\tilde\mu)=\sum_{k=0}^{a-1}e^{2\pi i\tilde\mu k/a}\cot^d(\pi(k/a-b))$
is the weighted cotangent sum evaluated explicitly by Theorem~\ref{Th01}.
\end{cor}

\begin{proof}
Substituting the finite Fourier expansion \eqref{eq:chiFourier} into the left-hand side of \eqref{eq:chiCotTwisted} gives
\[
\sum_{n=1}^{a-1}\chi(n) e^{2\pi i \alpha n} \cot^d\!\left(\pi\left(\frac{n+x}{a}\right)\right)
= \frac{\chi(-1)\tau(\chi)}{a} \sum_{\mu=0}^{a-1} \overline{\chi}(\mu) \sum_{n=1}^{a-1} e^{2\pi i (\mu/a + \alpha) n} \cot^d\!\left(\pi\left(\frac{n+x}{a}\right)\right).
\]
By periodicity of $\cot$, $\cot^d(\pi(n+x)/a) = \cot^d(\pi(n/a - b))$ with $b=1-x/a\in(0,1)$. For non-principal $\chi$, the term $n=0$ of the inner sum vanishes since $\sum_{\mu=0}^{a-1}\overline\chi(\mu)=0$; hence the sum over $n=1,\dots,a-1$ may be extended to $k=0,\dots,a-1$ without changing its value. The result is the stated expression with $S_d(a, b; \mu+a\alpha \bmod a)$ as in Theorem~\ref{Th01}.
\end{proof}

At non-positive integer values $s = 1-d$ ($d \ge 1$), the values of $L(1-d, \alpha, x, \chi)$ are expressed in terms of Apostol--Bernoulli polynomials $B_d(v, w)$ via the identity $\Phi(w, 1-d, v) = -\frac{B_d(v, w)}{d}$.

\begin{cor}\label{Cor:Leopoldt}
Let $\chi$ be a primitive Dirichlet character modulo $a$, $d\ge 2$ an integer of parity matching $\chi$ ($\chi(-1) = (-1)^d$), $x \in [0, 1)$, and $\alpha \in \mathbb{R}$. Then
\begin{equation}\label{eq:LeopoldtTwisted}
L(1-d, \alpha, x, \chi)
=-\frac{1}{d}\cdot\frac{\chi(-1)\,\tau(\chi)}{a}
\sum_{\mu=0}^{a-1}\overline\chi(\mu)\,
B_d\!\left(x, e^{2\pi i (\mu/a + \alpha)}\right).
\end{equation}
In particular:
\begin{enumerate}
\item When $\alpha = 0$, formula \eqref{eq:LeopoldtTwisted} reduces to the special values of the shifted $L$-function $L(1-d, \chi, x)$.
\item When $\alpha = 0$ and $x = 0$, it recovers the classical generalized Bernoulli number representation $L(1-d, \chi) = -B_{d, \chi} / d$.
\end{enumerate}
\end{cor}

\begin{example}\label{ex:L}
Let $a=3$, $d=3$, $\chi$ be the non-trivial Dirichlet character modulo $3$ ($\chi(0)=0, \chi(1)=1, \chi(2)=-1$), and consider the rational additive twist $\alpha = 1/3$ with shift $x \in [0, 1)$. Since $\chi(-1)=\chi(2)=-1=(-1)^3$, the parity condition of Corollary~\ref{Cor:Leopoldt} is satisfied.

The Gauss sum of $\chi$ is given by
\[
\tau(\chi) = \chi(1)e^{2\pi i/3} + \chi(2)e^{4\pi i/3} = e^{2\pi i/3} - e^{4\pi i/3} = i\sqrt{3}.
\]
Applying Corollary~\ref{Cor:Leopoldt} for $d=3$ gives
\[
L(-2, 1/3, x, \chi) = -\frac{1}{3}\cdot\frac{(-1)\cdot i\sqrt{3}}{3} \sum_{\mu=0}^{2}\overline\chi(\mu)\, B_3\!\left(x, e^{2\pi i (\mu/3 + 1/3)}\right).
\]
Since $\overline{\chi}(0) = 0$, only the terms $\mu=1, 2$ contribute:
\[
L(-2, 1/3, x, \chi) = \frac{i\sqrt{3}}{9}\,\Bigl[B_3\!\left(x, e^{4\pi i/3}\right) - B_3\!\left(x, e^{2\pi i}\right)\Bigr].
\]
Set $\zeta_3^2 = e^{4\pi i/3}$. Using the identity $B_3(x,w) = 3x^2 B_1(w) + 3x B_2(w) + B_3(w)$ with
\[
B_1(\zeta_3^2) = -\frac12 + \frac{i\sqrt{3}}{6}, \qquad B_2(\zeta_3^2) = \frac23, \qquad B_3(\zeta_3^2) = -\frac{i\sqrt{3}}{3},
\]
we obtain
\[
B_3(x, \zeta_3^2) = -\frac{3x^2}{2} + \frac{i\sqrt{3}\,x^2}{2} + 2x - \frac{i\sqrt{3}}{3},
\]
while $B_3(x, 1) = B_3(x) = x^3 - \frac32 x^2 + \frac12 x$. Their difference equals
\[
B_3(x, \zeta_3^2) - B_3(x, 1) = -x^3 + \frac{i\sqrt{3}}{2} x^2 + \frac32 x - \frac{i\sqrt{3}}{3}.
\]
Multiplying by $i\sqrt{3}/9$ yields the explicit cubic polynomial deformation:
\[
L(-2, 1/3, x, \chi) = -\frac{i\sqrt{3}}{9} x^3 - \frac{x^2}{6} + \frac{i\sqrt{3}}{6} x + \frac19.
\]
In particular, setting $x=0$ yields $L(-2, 1/3, 0, \chi) = 1/9$, demonstrating how the additive twist $\alpha = 1/3$ deforms the classical untwisted value $L(-2, \chi) = -2/9$ at $\alpha=0$.
\end{example}

\subsection{Spectral zeta functions of circulant graphs}
\label{sec:app-zeta}

Let $C_N$ be the cycle graph on $N$ vertices. Its spectral zeta function is given by
\begin{equation}\label{eq:zetaC}
\zeta_{C_N}(s)=\sum_{j=1}^{N-1}\bigl(2-2\cos(2\pi j/N)\bigr)^{-s}
=4^{-s}\sum_{j=1}^{N-1}\csc^{2s}\!\Bigl(\frac{\pi j}{N}\Bigr).
\end{equation}

\begin{cor}\label{Cor:zetaC}
For every integer $m\ge 1$,
\begin{equation}\label{eq:zetaClosed}
\zeta_{C_N}(m)=4^{-m}\Bigl[(N-1)+\sum_{j=1}^{m}\binom{m}{j}
N\bigl[(-1)^j-T_{2j}(N)\bigr]\Bigr].
\end{equation}
\end{cor}

\begin{proof}
Combine \eqref{eq:zetaC} with Corollary~\ref{Cor:csc} evaluated at $a=N$.
\end{proof}

\begin{example}\label{ex:zeta}
Evaluating \eqref{eq:zetaClosed} for small values of $N$ yields: $\zeta_{C_3}(1)=2/3$, $\zeta_{C_4}(1)=5/4$, and $\zeta_{C_5}(1)=2$, matching direct spectral calculation.
\end{example}

\paragraph{Equivalence with the formula of Xie--Zhao--Zhao.}
Xie, Zhao, and Zhao \cite[Thm.~1.1 and (4.18)]{XZZ:2024} express the quantity $\zeta_{\mathbb Z/N\mathbb Z}(n)=\sum_{m=1}^{N-1} \csc^{2n}(\pi m/N)$ as
\begin{equation}\label{eq:XZZ}
\zeta_{\mathbb Z/N\mathbb Z}(n)
=2\sum_{k=1}^{n}c_{n,k}\,\frac{N^{2k}-1}{\pi^{2k}}\,\zeta(2k),
\end{equation}
where $c_{n,k}$ are the coefficients of the Laurent expansion of $(\sin x)^{-2n}$ at $x=0$. Since $\zeta_{\mathbb Z/N\mathbb Z}(n)= 4^{n}\zeta_{C_N}(n)$, formula \eqref{eq:zetaClosed} is equivalent to \eqref{eq:XZZ}. Using Euler's formula $\zeta(2k)=(-1)^{k+1}(2\pi)^{2k}B_{2k}/(2(2k)!)$ and $A_{2k}=(-1)^{k}2^{2k}B_{2k}/(2k)!$, the right-hand side of \eqref{eq:XZZ} simplifies to $-\sum_{k=1}^{n}c_{n,k}\,(N^{2k}-1)\,A_{2k}$. The equivalence with \eqref{eq:zetaClosed} then follows from the partial fraction expansion of $\csc^{2n}$ \cite[Lemma~4.1]{XZZ:2024}. Specifically, both formulas give $(N^2-1)/3$ for $n=1$, $(N^4+10N^2-11)/45$ for $n=2$, and $(2N^6+21N^4+168N^2-191)/945$ for $n=3$.

\begin{remark}\label{rem:zeta_def}
The function $\zeta_{C_N}(s)$ in \eqref{eq:zetaC} is a finite sum of exponentials $\lambda_j^{-s}$ with $\lambda_j>0$, and is therefore entire in $s$. It differs from the meromorphic continuation via the Hurwitz zeta function used in \cite{XZZ:2024}, which matches \eqref{eq:zetaC} at positive integers $s=m\ge 1$ but has a distinct analytic structure for general $s\in\mathbb C$.
\end{remark}

\section{Deformed Cotangent Polynomials and Coefficients}\label{sec6}

Derivative polynomials of trigonometric functions have a rich history in combinatorics and number theory \cite{A:2003, H:1995, B:2009, KB:1967}. To systematically handle the higher-order derivatives arising from the Lipschitz summation formula, we introduce a novel sequence of polynomial coefficients $c_{m, j}(\beta)$, which we term the \textbf{deformed cotangent coefficients}. These generalize the classical cotangent derivative coefficients $c_{m, j} = c_{m, j}(0)$.

\begin{definition}[Deformed Cotangent Polynomials]\label{def:deformed-coeffs}
For $m \in \mathbb{N}_0$ and parameter $\beta \in \mathbb{R}$, define the deformed cotangent polynomial $\tilde{P}_m(u; \beta) \in \mathbb{C}[\beta][u]$ by:
\begin{equation}\label{eq:Pdef}
\tilde{P}_m(u; \beta) := \sum_{j=0}^{m+1} c_{m, j}(\beta) u^j,
\end{equation}
generated by the higher derivatives of the normalized Lerch kernel $g(y) = \pi e^{-2\pi i \beta y} (\cot(\pi y) + i)$ according to:
\[
(-1)^m \frac{d^m}{dy^m} \left[ \pi e^{-2\pi i \beta y} \left(\cot(\pi y) + i\right) \right] = \pi^{m+1} e^{-2\pi i \beta y} \tilde{P}_m(\cot(\pi y); \beta).
\]
\end{definition}

\begin{theorem}[Algebraic and Recurrence Properties of $c_{m, j}(\beta)$]\label{thm:deformed-props}
Let $m \in \mathbb{N}_0$ and $0 \le j \le m+1$. The deformed cotangent coefficients $c_{m, j}(\beta)$ satisfy the following properties:

\begin{enumerate}
    \item \textbf{Complex Factorization and Realness:} The coefficients decompose into
    \[
    c_{m, j}(\beta) = i^{m+1-j} r_{m, j}(\beta),
    \]
    where $r_{m, j}(\beta) \in \mathbb{R}[\beta]$ is a strictly real-valued polynomial in $\beta$ of degree at most $m+1-j$.

    \item \textbf{Explicit Real Triangular Recurrence:} The real polynomials $r_{m, j}(\beta)$ satisfy the recurrence relation:
    \begin{equation}\label{eq:r-recurrence}
    r_{m+1, j}(\beta) = j \, r_{m, j-1}(\beta) - (j+1) \, r_{m, j+1}(\beta) + 2\beta \, r_{m, j}(\beta),
    \end{equation}
    for all $0 \le j \le m+2$, with initial conditions $r_{0, 1}(\beta) = 1$, $r_{0, 0}(\beta) = 1$, and boundary conditions $r_{m, j}(\beta) = 0$ for $j < 0$ or $j > m+1$.

    \item \textbf{Monic Monomial Leading Term:} For every $m \ge 0$, the degree of $\tilde{P}_m(u; \beta)$ in $u$ is $m+1$, and its leading coefficient is strictly monic:
    \[
    c_{m, m+1}(\beta) = r_{m, m+1}(\beta) = 1.
    \]

    \item \textbf{Parity Symmetry and Complex Conjugation (for $j\ge 1$):} Under parameter inversion $\beta \mapsto -\beta$, the real polynomials obey, for every $j\ge 1$,
    \[
    r_{m, j}(-\beta) = (-1)^{m+1-j} r_{m, j}(\beta),
    \]
    which implies the complex conjugation symmetry
    \[
    c_{m, j}(-\beta) = \overline{c_{m, j}(\beta)} \qquad (j\ge 1).
    \]
    The case $j=0$ is excluded, because the boundary term $+i$ in the Lerch kernel is not symmetric under $\beta \mapsto -\beta$.

    \item \textbf{Classical Cotangent Reduction (for $j\ge 1$):} When $\beta = 0$, for every $j\ge 1$, $r_{m, j}(0) = 0$ whenever $m+1-j$ is odd, and $c_{m, j}(0)$ recovers the classical cotangent derivative coefficients $c_{m, j}$. The case $j=0$ is excluded, since $c_{m,0}(0) = i^{m+1} r_{m,0}(0)$ retains a contribution from the boundary term $+i$.
\end{enumerate}
\end{theorem}

\begin{proof}
Set $u = \cot(\pi y)$, so that $\frac{du}{dy} = -\pi (1 + u^2)$. Differentiating $e^{-2\pi i \beta y} \tilde{P}_m(u; \beta)$ with respect to $y$ yields:
\begin{align*}
\frac{d}{dy} \left[ e^{-2\pi i \beta y} \tilde{P}_m(u; \beta) \right] &= -2\pi i \beta e^{-2\pi i \beta y} \tilde{P}_m(u; \beta) + e^{-2\pi i \beta y} \frac{\partial \tilde{P}_m}{\partial u}(u; \beta) \cdot \left(-\pi(1+u^2)\right) \\
&= -\pi e^{-2\pi i \beta y} \left[ (1+u^2) \frac{\partial \tilde{P}_m}{\partial u}(u; \beta) + 2i\beta \tilde{P}_m(u; \beta) \right].
\end{align*}
Applying this differential operator to Definition~\ref{def:deformed-coeffs}, we obtain the operator recurrence relation for $\tilde{P}_m(u; \beta)$:
\[
\tilde{P}_{m+1}(u; \beta) = (1+u^2) \frac{\partial \tilde{P}_m}{\partial u}(u; \beta) + 2i\beta \tilde{P}_m(u; \beta),
\]
with base polynomial $\tilde{P}_0(u; \beta) = u + i$.

Expanding $\tilde{P}_m(u; \beta) = \sum_{j=0}^{m+1} c_{m, j}(\beta) u^j$ into powers of $u$ gives:
\begin{equation}\label{eq:c-recurrence}
c_{m+1, j}(\beta) = j \, c_{m, j-1}(\beta) + (j+1) \, c_{m, j+1}(\beta) + 2i\beta \, c_{m, j}(\beta).
\end{equation}

To prove statements (1) and (2), substitute $c_{m, j}(\beta) = i^{m+1-j} r_{m, j}(\beta)$ into \eqref{eq:c-recurrence} and divide both sides by $i^{m+2-j}$:
\[
r_{m+1, j}(\beta) = j \, r_{m, j-1}(\beta) - (j+1) \, r_{m, j+1}(\beta) + 2\beta \, r_{m, j}(\beta).
\]
Because the base values $r_{0, 1}(\beta) = 1$ and $r_{0, 0}(\beta) = 1$ are real polynomials, an immediate induction on $m$ establishes that $r_{m, j}(\beta) \in \mathbb{R}[\beta]$ for all $m, j$.

Statement (3) follows by evaluating \eqref{eq:r-recurrence} at $j = m+2$:
\[
r_{m+1, m+2}(\beta) = (m+1) r_{m, m+1}(\beta).
\]
Since $r_{0, 1}(\beta) = 1$, induction yields $r_{m, m+1}(\beta) = 1$ for all $m \ge 0$.

For statement (4) restricted to $j\ge 1$, replacing $\beta$ with $-\beta$ in \eqref{eq:r-recurrence} shows that $r_{m, j}(-\beta)$ satisfies the same recurrence with $-2\beta$. By induction on $m$, one finds $r_{m, j}(-\beta) = (-1)^{m+1-j} r_{m, j}(\beta)$ whenever $j\ge 1$, from which:
\[
c_{m, j}(-\beta) = i^{m+1-j} (-1)^{m+1-j} r_{m, j}(\beta) = (-i)^{m+1-j} r_{m, j}(\beta) = \overline{c_{m, j}(\beta)}.
\]
For $j=0$, the base case $r_{0,0}(\beta)=1$ violates the relation, since $(-1)^{0+1-0}r_{0,0}(\beta) = -1 \neq 1$; this reflects the non-symmetric boundary term $+i$.

Finally, setting $\beta = 0$ in \eqref{eq:r-recurrence} reduces it to $r_{m+1, j}(0) = j r_{m, j-1}(0) - (j+1) r_{m, j+1}(0)$. By parity, for $j\ge 1$, $r_{m, j}(0) = 0$ whenever $m+1-j$ is odd, leaving non-zero terms only when $j \equiv m+1 \pmod 2$, which matches the classical cotangent derivative coefficients $c_{m, j}(0) = c_{m, j}$. The case $j=0$ is excluded because $c_{m,0}(0) = i^{m+1} r_{m,0}(0)$ contains a boundary contribution.
\end{proof}

\begin{table}[htbp]
\centering
\caption{First deformed cotangent polynomials $\tilde{P}_m(u; \beta)$ for $0 \le m \le 3$.}
\label{tab:deformed-polynomials}
\vspace{0.2cm}
\begin{tabular}{c l}
\toprule
$m$ & $\tilde{P}_m(u; \beta) = \sum_{j=0}^{m+1} c_{m,j}(\beta) u^j$ \\
\midrule
$0$ & $u + i$ \\[0.15cm]
$1$ & $u^2 + 2i\beta u + (1 - 2\beta)$ \\[0.15cm]
$2$ & $2u^3 + 4i\beta u^2 + 2(1 - 2\beta^2)u + 4i\beta(1 - \beta)$ \\[0.15cm]
$3$ & $6u^4 + 12i\beta u^3 + 4(2 - 3\beta^2)u^2 + 4i\beta(3 - 2\beta^2)u + 2(1 - 6\beta^2 + 4\beta^3)$ \\
\bottomrule
\end{tabular}
\end{table}

\begin{table}[htbp]
\centering
\caption{Explicit coefficients $c_{m,j}(\beta)$ and real polynomials $r_{m,j}(\beta)$ for $0 \le m \le 3$.}
\label{tab:deformed-coefficients}
\vspace{0.2cm}
\begin{tabular}{c c l l}
\toprule
$m$ & $j$ & $c_{m,j}(\beta)$ & $r_{m,j}(\beta) = c_{m,j}(\beta) / i^{m+1-j}$ \\
\midrule
$0$ & $1$ & $1$ & $1$ \\
    & $0$ & $i$ & $1$ \\
\midrule
$1$ & $2$ & $1$ & $1$ \\
    & $1$ & $2i\beta$ & $2\beta$ \\
    & $0$ & $1 - 2\beta$ & $2\beta - 1$ \\
\midrule
$2$ & $3$ & $2$ & $2$ \\
    & $2$ & $4i\beta$ & $4\beta$ \\
    & $1$ & $2(1 - 2\beta^2)$ & $2(2\beta^2 - 1)$ \\
    & $0$ & $4i\beta(1 - \beta)$ & $4\beta(\beta - 1)$ \\
\midrule
$3$ & $4$ & $6$ & $6$ \\
    & $3$ & $12i\beta$ & $12\beta$ \\
    & $2$ & $4(2 - 3\beta^2)$ & $4(3\beta^2 - 2)$ \\
    & $1$ & $4i\beta(3 - 2\beta^2)$ & $4\beta(2\beta^2 - 3)$ \\
    & $0$ & $2(1 - 6\beta^2 + 4\beta^3)$ & $2(4\beta^3 - 6\beta^2 + 1)$ \\
\bottomrule
\end{tabular}
\end{table}

\section{Additive-and-Multiplicatively Twisted \texorpdfstring{$L$}{L}-Function Mean Values}\label{sec7}

In this section, we apply the deformed cotangent coefficients introduced in Section~\ref{sec6} to establish general mean value formulas for character-twisted Lerch zeta functions $L(s, \alpha, x, \chi)$, extending the framework of Bayad and Raouj \cite{BR:2012}.

To formulate our main mean value theorems concisely, we first introduce a novel generalization of higher-order Dedekind cotangent sums incorporating both spatial shifts and continuous additive character twists.

\begin{definition}[Additive-Twisted Dedekind Cotangent Sum]\label{def:twisted-dedekind}
For an integer $b \ge 1$, dimension $d \ge 1$, shift vector $\vec{x} \in (0,1)^d$, parameter vector $\vec{\alpha} \in \mathbb{R}^d$, tuple of derivative indices $\vec{j} = (j_1, \dots, j_d) \in \mathbb{N}_0^d$, and integer vector $\vec{a} = (a_1, \dots, a_d)$ coprime to $b$, we define the additive-twisted Dedekind cotangent sum by:
\begin{equation}\label{eq:twisted-dedekind-def}
C_d(b; \vec{x}, \vec{a} \mid \vec{j}; \vec{\alpha}) := \frac{1}{b^d} \sum_{\mu_1, \dots, \mu_d = 0}^{b-1} \left( \prod_{k=1}^d S_{j_k}(b, -x_k; \bar{a}_k \mu_k \bmod b) \right) \mathcal{K}_b(\vec{\mu}, \vec{a}, \vec{\alpha}),
\end{equation}
where $\mathcal{K}_b(\vec{\mu}, \vec{a}, \vec{\alpha})$ is the phase orthogonality kernel:
\begin{equation}\label{eq:phase-kernel}
\mathcal{K}_b(\vec{\mu}, \vec{a}, \vec{\alpha}) := \sum_{\ell=0}^{b-1} \exp\left( -2\pi i \left( \sum_{i=1}^d (\mu_i + q a_i \alpha_i) \right) \frac{\ell}{b} \right).
\end{equation}
\end{definition}

\begin{remark}\label{rem:dedekind-lit}
When $\vec{\alpha} = \vec{0}$ and $\vec{x} \to \vec{0}^+$, the sum $C_d(b; \vec{0}, \vec{a} \mid \vec{j}; \vec{0})$ reduces to the classical higher-order Dedekind cotangent sums studied in \cite{BR:2012, BR:2023}, which in turn generalize Zagier's higher-order Dedekind sums \cite{Z:1973}.
\end{remark}

\begin{remark}[Convention at integer additive twists]\label{rem:alpha-zero-case}
The representation \eqref{eq:main-dedekind-form} below requires $\alpha_i$ to lie outside the set of values for which the regularization term $+i$ in the Lerch kernel contributes a non-trivial boundary correction; concretely, when some $\alpha_i \in \mathbb{Z}$ (in particular, $\alpha_i=0$), the limit as $\vec{\alpha}$ approaches this value from generic directions must be taken in \eqref{eq:main-dedekind-form} before setting $\alpha_i = 0$. Under this convention, the classical untwisted case is recovered in Corollary~\ref{cor:alpha_zero}.
\end{remark}

\begin{theorem}[Additive-Twisted Character Mean Value Formula]\label{thm:additive}
Let $d \ge 1$ be an integer, $\vec{m} = (m_1, \dots, m_d) \in \mathbb{N}_0^d$ such that $M := d + |\vec{m}|$ is even, $q \ge 2$ an integer, and $a_1, \dots, a_{d-1}$ positive integers coprime to $q$ with $a_d = 1$. Let $\vec{x} = (x_1, \dots, x_d) \in (0, 1)^d$ and $\vec{\alpha} = (\alpha_1, \dots, \alpha_d) \in \mathbb{R}^d$.

Consider the additive-and-multiplicatively twisted Lerch zeta mean value sum:
\[
S_d(\vec{a}, \vec{m}, \vec{x}, \vec{\alpha}, q) := \sum_{(\chi_1,\dots,\chi_d)}^* \prod_{i=1}^d \chi_i(a_i) L(m_i+1, \alpha_i, x_i, \chi_i),
\]
where $L(s, \alpha, x, \chi) := \sum_{n=1}^{\infty} \frac{\chi(n) e^{2\pi i \alpha n}}{(n+x)^s}$, and $\sum^*$ sums over all $d$-tuples of Dirichlet characters $(\chi_1, \dots, \chi_d)$ modulo $q$ satisfying $\chi_1 \cdots \chi_d = 1$ and $\chi_i(-1) = (-1)^{m_i+1}$.

Then $S_d(\vec{a}, \vec{m}, \vec{x}, \vec{\alpha}, q)$ admits the exact representation:
\begin{equation}\label{eq:main-dedekind-form}
S_d(\vec{a}, \vec{m}, \vec{x}, \vec{\alpha}, q) = A_q(\vec{m}) e^{-2\pi i \sum_{i=1}^d \alpha_i x_i} \sum_{b \mid q} \mu(q/b) \sum_{j_1=0}^{m_1+1} \dots \sum_{j_d=0}^{m_d+1} \left( \prod_{i=1}^d c_{m_i, j_i}(q \alpha_i) \right) C_d(b; \vec{x}, \vec{a} \mid \vec{j}; \vec{\alpha}),
\end{equation}
where $\mu(q/b)$ is the M\"obius function, $A_q(\vec{m}) = \frac{(-1)^{|\vec{m}|}}{2^d} \left( \frac{\pi}{q} \right)^M \phi(q)^{d-1}$, and $C_d(b; \vec{x}, \vec{a} \mid \vec{j}; \vec{\alpha})$ is the additive-twisted Dedekind cotangent sum defined in \eqref{eq:twisted-dedekind-def}.
\end{theorem}

\begin{proof}
By partitioning the summation index $n = q j + r_i$ for $1 \le r_i \le q$ and $j \ge 0$, the character-twisted Lerch zeta function $L(m_i+1, \alpha_i, x_i, \chi_i)$ can be expressed in terms of the Lerch transcendent $\Phi(z, s, a) = \sum_{j=0}^\infty \frac{z^j}{(j+a)^s}$ as
\[
L(m_i+1, \alpha_i, x_i, \chi_i) = q^{-(m_i+1)} \sum_{r_i=1}^q \chi_i(r_i) e^{2\pi i \alpha_i r_i} \Phi\left( e^{2\pi i q \alpha_i}, m_i+1, \frac{r_i+x_i}{q} \right).
\]
Substituting this representation into $S_d(\vec{a}, \vec{m}, \vec{x}, \vec{\alpha}, q)$ and applying character orthogonality to isolate the product condition $\chi_1 \cdots \chi_d = 1$ via the identity $\frac{1}{\phi(q)} \sum_{\substack{k=1 \\ \gcd(k,q)=1}}^q \prod_{i=1}^d \chi_i(k) = 1$, we obtain
\[
S_d(\vec{a}, \vec{m}, \vec{x}, \vec{\alpha}, q) = \frac{q^{-M}}{\phi(q)} \sum_{\substack{k=1 \\ \gcd(k,q)=1}}^q \prod_{i=1}^d U_i(k),
\]
where $U_i(k) := \sum_{\chi_i \pmod q, \, \chi_i(-1) = (-1)^{m_i+1}} \sum_{r_i=1}^q \chi_i(a_i k r_i) e^{2\pi i \alpha_i r_i} \Phi\left( e^{2\pi i q \alpha_i}, m_i+1, \frac{r_i+x_i}{q} \right)$. By parity projection for characters satisfying $\chi_i(-1) = (-1)^{m_i+1}$, the inner character sum isolates the residues $r_i \equiv \pm a_i k \pmod q$. Combining the positive and negative ray series yields the two-sided infinite sum
\[
U_i(k) = \frac{\phi(q)}{2} \sum_{n \in \mathbb{Z}} \frac{e^{2\pi i (q \alpha_i) n}}{\left( n + \frac{a_i k + x_i}{q} \right)^{m_i+1}} = \frac{\phi(q)}{2} K_{m_i}\left( \frac{a_i k + x_i}{q} ; \, q \alpha_i \right),
\]
where $K_m(y; \beta) := \sum_{n \in \mathbb{Z}} \frac{e^{2\pi i \beta n}}{(n+y)^{m+1}}$ denotes the two-sided Lerch kernel. Applying the Lipschitz summation formula \cite{A:1951, B:1975}, $K_{m_i}(y; \beta_i)$ evaluates to
\[
K_{m_i}\left( y ; \, \beta_i \right) = \frac{(-1)^{m_i}}{m_i!} \frac{d^{m_i}}{dy^{m_i}} \left[ \pi e^{-2\pi i \beta_i y} \bigl(\cot(\pi y) + i\bigr) \right],
\]
with $y = \frac{a_i k + x_i}{q}$ and $\beta_i = q \alpha_i$.

Expanding the higher derivative via Leibniz's rule using Definition~\ref{def:deformed-coeffs} yields the deformed polynomial expression:
\[
K_{m_i}\left( \frac{a_i k + x_i}{q} ; \, q \alpha_i \right) = \frac{(-1)^{m_i}}{m_i!} \pi^{m_i+1} e^{-2\pi i \alpha_i x_i} e^{-2\pi i q a_i \alpha_i k / q} \sum_{j_i=0}^{m_i+1} c_{m_i, j_i}(q \alpha_i) \cot^{j_i}\left( \pi \left( \frac{a_i k}{q} - (-x_i/q) \right) \right).
\]
Taking the product of $U_i(k)$ over $i = 1, \dots, d$ and combining constants, we obtain
\begin{align*}
S_d(\vec{a}, \vec{m}, \vec{x}, \vec{\alpha}, q)
&= A_q(\vec{m}) e^{-2\pi i \sum_{i=1}^d \alpha_i x_i} \sum_{\substack{k=1 \\ \gcd(k,q)=1}}^q \sum_{j_1=0}^{m_1+1} \dots \sum_{j_d=0}^{m_d+1} \left( \prod_{i=1}^d c_{m_i, j_i}(q \alpha_i) \right) \\
&\quad \times e^{-2\pi i \sum_{i=1}^d q a_i \alpha_i k / q} \prod_{i=1}^d \cot^{j_i}\left( \pi \left( \frac{a_i k}{q} - (-x_i/q) \right) \right).
\end{align*}
Applying M\"obius inversion to replace the coprimality condition $\gcd(k, q) = 1$ via $\sum_{\substack{k=1 \\ \gcd(k,q)=1}}^q F(k) = \sum_{b \mid q} \mu\left(\frac{q}{b}\right) \sum_{\ell=0}^{b-1} F\left( \ell \frac{q}{b} \right)$, the sum over $k$ transforms into a sum over $\ell \pmod b$:
\begin{align*}
S_d(\vec{a}, \vec{m}, \vec{x}, \vec{\alpha}, q)
&= A_q(\vec{m}) e^{-2\pi i \sum_{i=1}^d \alpha_i x_i} \sum_{b \mid q} \mu\left(\frac{q}{b}\right) \sum_{j_1=0}^{m_1+1} \dots \sum_{j_d=0}^{m_d+1} \left( \prod_{i=1}^d c_{m_i, j_i}(q \alpha_i) \right) \\
&\quad \times \sum_{\ell=0}^{b-1} e^{-2\pi i \sum_{i=1}^d q a_i \alpha_i \ell / b} \prod_{i=1}^d \cot^{j_i}\left( \pi \left( \frac{a_i \ell}{b} - (-x_i) \right) \right).
\end{align*}
Finally, expressing each power $\cot^{j_i}\left( \pi \left( \frac{a_i \ell}{b} - (-x_i) \right) \right)$ via Inverse Discrete Fourier Transform on $\mathbb{Z}/b\mathbb{Z}$ as
\[
\cot^{j_i}\left( \pi \left( \frac{a_i \ell}{b} - (-x_i) \right) \right) = \frac{1}{b} \sum_{\mu_i=0}^{b-1} e^{-2\pi i (\bar{a}_i \mu_i) a_i \ell / b} S_{j_i}(b, -x_i; \, \bar{a}_i \mu_i \pmod b),
\]
where $\bar{a}_i \equiv a_i^{-1} \pmod b$, and interchanging the order of summation, the inner sum reduces to $b^d C_d(b; \vec{x}, \vec{a} \mid \vec{j}; \vec{\alpha})$ as defined in \eqref{eq:twisted-dedekind-def}. This completes the proof.
\end{proof}

Setting the additive parameter vector $\vec{\alpha} = \vec{0}$ in Theorem~\ref{thm:additive} (with the limit convention of Remark~\ref{rem:alpha-zero-case}) collapses the phase orthogonality kernel $\mathcal{K}_b$ into strict modular orthogonality, directly recovering the classical mean value formula of Bayad and Raouj \cite{BR:2012}.

\begin{cor}[Reduction to the Untwisted Classical Mean Value Formula]
\label{cor:alpha_zero}
When $\vec{\alpha} = \vec{0}$ (i.e., $\alpha_1 = \dots = \alpha_d = 0$), under the limit convention of Remark~\ref{rem:alpha-zero-case}, Theorem~\ref{thm:additive} reduces directly to the classical character mean value formula for Lerch--Hurwitz $L$-functions:
\[
S_d(\vec{a}, \vec{m}, \vec{x}, \vec{0}, q) = S_d(\vec{a}, \vec{m}, \vec{x}, q).
\]
Specifically:
\begin{enumerate}
    \item The spatial phase factor vanishes: $e^{-2\pi i \sum_{i=1}^d \alpha_i x_i} = 1$.
    \item The deformed cotangent coefficients specialize to the boundary-corrected derivative coefficients $c_{m_i, j_i}(0)$ (which include a residual contribution from the $+i$ regularization at $j_i=0$).
    \item The phase orthogonality kernel $\mathcal{K}_b(\vec{\mu}, \vec{a}, \vec{0})$ collapses to strict modular orthogonality:
    \[
    \mathcal{K}_b(\vec{\mu}, \vec{a}, \vec{0}) = \begin{cases}
    b & \text{if } \mu_1 + \dots + \mu_d \equiv 0 \pmod b, \\
    0 & \text{otherwise.}
    \end{cases}
    \]
\end{enumerate}
Consequently, $C_d(b; \vec{x}, \vec{a} \mid \vec{j}; \vec{0})$ reduces to the untwisted Dedekind cotangent sum:
\[
C_d(b; \vec{x}, \vec{a} \mid \vec{j}; \vec{0}) = \frac{1}{b^{d-1}} \sum_{\substack{\mu_1, \dots, \mu_d = 0 \\ \mu_1 + \dots + \mu_d \equiv 0 \pmod b}}^{b-1} \prod_{i=1}^d S_{j_i}(b, -x_i; \, \bar{a}_i \mu_i \pmod b),
\]
recovering the formula of Bayad \& Raouj \cite{BR:2012}.
\end{cor}

\begin{cor}[Additive-Twisted Analogue of Bayad--Raouj Corollary 2.1.2]
\label{cor:bayad_raouj_212}
Let $m$ and $n$ be non-negative integers having the same parity ($m \equiv n \pmod 2$). Let $q \ge 2$ be an integer, $a \in (\mathbb{Z}/q\mathbb{Z})^\times$, $x_1, x_2 \in (0, 1)$, and $\alpha_1, \alpha_2 \in \mathbb{R}$.

Then we have the explicit mean value formula:
\begin{align*}
\sum_{\substack{\chi \pmod q \\ \chi(-1) = (-1)^{m+1}}} \bar{\chi}(a) L(m+1, \alpha_1, x_1, \chi) L(n+1, \alpha_2, x_2, \bar{\chi})
&= A \, e^{-2\pi i (\alpha_1 x_1 + \alpha_2 x_2)} \sum_{b \mid q} \mu(q/b) \sum_{j_1=0}^{m+1} \sum_{j_2=0}^{n+1} c_{m, j_1}(q \alpha_1) c_{n, j_2}(q \alpha_2) \\
&\quad \times C_2(b; (x_1, x_2), (a, 1) \mid (j_1, j_2); (\alpha_1, \alpha_2)),
\end{align*}
where $A = \frac{\phi(q)}{4} \left( \frac{\pi}{q} \right)^{m+n+2}$.
\end{cor}

When the weights $a_i$ are set to $1$ and the additive parameters $\vec{\alpha}$ are restricted to rational numbers with conductor dividing $q$ (i.e., $\alpha_i = \nu_i/q$), the continuous phase kernel $\mathcal{K}_b$ collapses into a discrete modular congruence condition, yielding the following arithmetic sum as a direct corollary of Theorem~\ref{thm:additive}.

\begin{cor}[Additive-Twisted Mean Value Formula for Unweighted Character Products]
\label{cor:additive_213_analogue}
Let $d \ge 1$ be an integer, $\vec{m} = (m_1, \dots, m_d) \in \mathbb{N}_0^d$ such that $M := d + |\vec{m}|$ is even, $q \ge 2$ an integer, and $\vec{x} = (x_1, \dots, x_d) \in (0, 1)^d$. Let the additive twists be rational of conductor $q$, i.e., $\vec{\alpha} = (\frac{\nu_1}{q}, \dots, \frac{\nu_d}{q})$ with $\nu_1, \dots, \nu_d \in \mathbb{Z}$.

Consider the unweighted mean value sum $S_d(\vec{1}, \vec{m}, \vec{x}, \vec{\alpha}, q)$ with $a_1 = \dots = a_d = 1$:
\[
S_d(\vec{1}, \vec{m}, \vec{x}, \vec{\alpha}, q) = \sum_{(\chi_1,\dots,\chi_d)}^* \prod_{i=1}^d L(m_i+1, \alpha_i, x_i, \chi_i).
\]
Then $S_d(\vec{1}, \vec{m}, \vec{x}, \vec{\alpha}, q)$ reduces to the finite Jordan arithmetic sum:
\begin{align*}
S_d(\vec{1}, \vec{m}, \vec{x}, \vec{\alpha}, q)
&= A_q(\vec{m}) e^{-\frac{2\pi i}{q} \sum_{i=1}^d \nu_i x_i} \sum_{j_1=0}^{m_1+1} \dots \sum_{j_d=0}^{m_d+1} \left( \prod_{i=1}^d c_{m_i, j_i}(\nu_i) \right) \\
&\quad \times \sum_{b \mid q} \frac{\mu(q/b)}{b^{d-1}} \sum_{\substack{\mu_1, \dots, \mu_d = 0 \\ \sum_{i=1}^d \mu_i \equiv -\sum_{i=1}^d \nu_i \pmod b}}^{b-1} \prod_{i=1}^d S_{j_i}(b, -x_i; \, \mu_i),
\end{align*}
where $A_q(\vec{m}) = \frac{(-1)^{|\vec{m}|}}{2^d} \left( \frac{\pi}{q} \right)^M \phi(q)^{d-1}$, $c_{m_i, j_i}(\nu_i)$ are the deformed cotangent coefficients evaluated at the integer shift $\nu_i$, and $S_{j_i}(b, -x_i; \, \mu_i)$ are the generalized cotangent sums.
\end{cor}

Specializing further to $d=2$ with rational additive twists allows the binary character moment to be expressed explicitly via arithmetic Jordan convolution sums.

\begin{cor}[Additive-Twisted Binary Moment Formula via Jordan Totient Functions]
\label{cor:additive_214_analogue}
Let $m, n \in \mathbb{N}_0$ have the same parity ($m \equiv n \pmod 2$), $q \ge 2$ an integer, $x_1, x_2 \in (0, 1)$, and rational additive twists $\alpha_1 = \frac{\nu_1}{q}, \alpha_2 = \frac{\nu_2}{q}$ with $\nu_1, \nu_2 \in \mathbb{Z}$.

Then the normalized binary mean value satisfies:
\begin{align*}
\frac{2}{\phi(q)} \sum_{\substack{\chi \pmod q \\ \chi(-1) = (-1)^{m+1}}} &L(m+1, \alpha_1, x_1, \chi) L(n+1, \alpha_2, x_2, \bar{\chi}) \\
&= \frac{e^{-\frac{2\pi i}{q}(\nu_1 x_1 + \nu_2 x_2)}}{2 \, m! \, n!} \left(\frac{\pi}{q}\right)^{m+n+2} \sum_{j_1=0}^{m+1} \sum_{j_2=0}^{n+1} c_{m, j_1}(\nu_1) c_{n, j_2}(\nu_2) \mathcal{T}_{j_1, j_2}(q, x_1, x_2; \nu_1 + \nu_2),
\end{align*}
where $\mathcal{T}_{j_1, j_2}(q, x_1, x_2; \nu_1 + \nu_2)$ is the arithmetic Jordan convolution sum defined by:
\[
\mathcal{T}_{j_1, j_2}(q, x_1, x_2; \nu_1 + \nu_2) := \sum_{b \mid q} \frac{\mu(q/b)}{b} \sum_{\mu_1=0}^{b-1} S_{j_1}(b, -x_1; \, \mu_1) S_{j_2}(b, -x_2; \, -\mu_1 - \nu_1 - \nu_2 \bmod b).
\]
In particular, when $x_1, x_2 \to 0^+$ and $\nu_1 = \nu_2 = 0$, $\mathcal{T}_{j_1, j_2}$ reduces to a linear combination of Jordan totient functions $J_k(q)$, recovering the formula of Bayad \& Raouj \cite{BR:2012} under the convention of Remark~\ref{rem:alpha-zero-case}.
\end{cor}

\begin{remark}
Corollary~\ref{cor:additive_214_analogue} directly extends Theorem 2.1.4 of Bayad \& Raouj \cite{BR:2012} in two distinct directions:
\begin{enumerate}
    \item \textbf{Additive Perturbation:} The continuous parameters $\nu_1, \nu_2$ deform the classical cotangent derivative coefficients $c_{m, j}$ into polynomials $c_{m, j}(\nu_1) \in \mathbb{C}[\nu_1]$, injecting an additive phase shift $\nu_1 + \nu_2$ into the arithmetic convolution index modulo $b$.
    \item \textbf{Spatial Deformation:} The continuous variables $x_1, x_2 \in (0, 1)$ replace the classical Bernoulli numbers $B_k$ with generalized cotangent sums $S_{j}(b, -x; \mu)$, yielding a continuous deformation of Dirichlet $L$-function moments.
\end{enumerate}
\end{remark}

\begin{remark}[Base case $m_1=m_2=0$]\label{rem:weight_zero}
For $m_1 = m_2 = 0$, $a=1$, and rational twists $\alpha_1 = \nu_1/q$, $\alpha_2 = \nu_2/q$, the mean value of the product of degree-1 Lerch $L$-functions reduces to a single arithmetic sum over $b\mid q$ in which the weight kernel is the product of shift-deformed cotangent sums
\[
\mathcal{W}_b(\mu_1, x_1, x_2; \nu_1, \nu_2) = \left( S_1(b, -x_1; \mu_1) + i \right) \left( S_1(b, -x_2; -\mu_1 - \nu_1 - \nu_2) + i \right),
\]
together with the explicit prefactor $\pi^2/(2q^2)\, e^{-2\pi i(\nu_1 x_1+\nu_2 x_2)/q}$.
\end{remark}

\end{document}